\documentclass{article}
\usepackage{amsmath,amsthm,amsfonts,amssymb,graphicx,hyperref}
\usepackage[shortlabels]{enumitem}
\usepackage{textcomp}

\begin{document}

\title{Bluffing in Scrabble}
\author{Nick Ballard \and Timothy Y. Chow}
\date{August 2026}

\maketitle

\section{The Puzzle}

\begin{figure}
\caption{You hold MKNOSYZ with DEFJLLLQW unseen. What's your play?}
\label{fig:puzzle}
\begin{center}
\includegraphics[scale=.7]{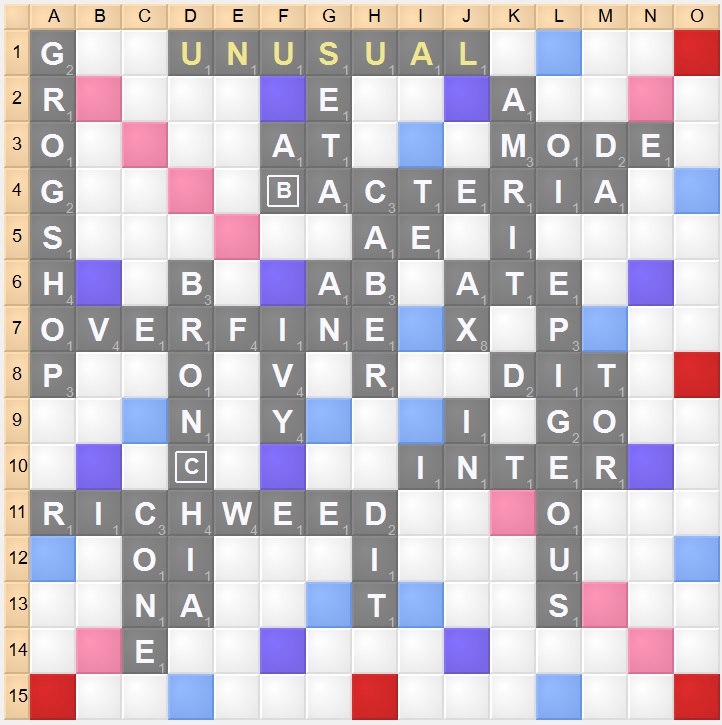}
\end{center}
\end{figure}

Imagine that you're playing Scrabble.
The game is nearing its end and you're in
dire straits. On your last turn, you bingoed with GROGSHOP,
but then your opponent immediately bingoed
back with UNUSUAL,
taking a commanding lead of 476 points to your 344 points.
See \autoref{fig:puzzle}.
On your rack, you have MKNOSYZ. The tiles unseen
to you are DEFJLLLQW, two of which are in the bag and
seven of which are on your opponent's rack.
What's your best play?

The surprising answer, subject to some assumptions that we explain below,
is that \emph{you should bluff} with probability~$2/3$.
We believe that \autoref{fig:puzzle} is the first published example
of a Scrabble position for which an equilibrium strategy
is necessarily a mixed strategy.
But before we explain in more detail what we mean,
let us briefly review some of the rules of Scrabble.

\section{Scrabble Scoring}

The rules of Scrabble\footnote{SCRABBLE\textsuperscript{\tiny\textregistered}
Brand Crossword Game is a registered trademark.
Throughout this paper, for brevity, we refer to the game
simply as ``Scrabble.''} are available from many sources,
but we should at least remind the reader of some key points.
The game begins with a publicly known set of tiles in an opaque bag.
To start the game, each of the two players draws seven random tiles,
which are placed on the player's rack
so that the opponent cannot see them
(the order in which the tiles are placed on the rack is immaterial).
Players take turns putting tiles from their rack onto the board,
scoring points by forming valid words,
where \emph{valid} means that the word
is listed in a pre-specified \emph{lexicon}.
The puzzle in \autoref{fig:puzzle} assumes the lexicon known as NWL2023,
which is the 2023 revision of the English word list
used in North American tournament play.

If you play all seven of your tiles at once,
it is called a \emph{bingo} and you score 50 extra points
on top of the points you earn from the word(s) you form.
Regardless of the number of tiles you play,
at the end of your turn, you draw fresh tiles randomly from the bag
so that you again have seven tiles on your rack,
unless there are too few tiles left in the bag,
in which case you draw as many tiles as you can.
If the bag is empty at the start of your turn, and you
play all the remaining tiles on your rack, then you
\emph{go out}, ending the game, and 
adding to your score twice the point value of any tiles
remaining on your opponent's rack.
The player with the higher score wins;
if both players have exactly the same score,
then the game ends in a tie.

In a typical Scrabble tournament, all the players play the same
number of games, earning one point for a win, zero points for a loss,
and half a point for a tie.
The player with the most points at the end of the tournament wins.
If two or more players have the same number of points,
then \emph{spread} is typically used as a tiebreaker.
\emph{Spread} refers to the amount by
which you outscored your opponents
(this amount is negative when you lose),
summed over the whole tournament.
However, in this paper, we ignore spread.

For notational purposes, the rows of a Scrabble board are
numbered 1--15, and the columns are labeled A--O.
The location of a play is specified by the coordinates
of the first letter of the word
(even if that letter is already on the board).
For horizontal words, the number is given first;
for vertical words, the letter is given first.

\section{Preliminary Analysis of the Puzzle}
\label{sec:preliminary}

Let us return to the puzzle.
You are certainly a big underdog and will probably lose,
so the only question is whether you have
any chances at all of saving the game,
assuming your opponent plays perfectly.
Given that you trail by 132 points,
your only hope of tying or winning is to bingo.
(For a justification of this claim,
see Appendix~\ref{app:analysis}.)
You currently have no playable bingo,
and if you play more than two tiles then
you will have at most six tiles on your rack next turn,
making a bingo impossible.
So, do you have any one-tile or two-tile plays
that might let you save the game with a bingo
on your next turn?

\begin{figure}[ht]
\caption{8K DITZ (left) Sets Up a Dream Sequence (right)}
\label{fig:dream}
\begin{center}
\includegraphics[scale=.45]{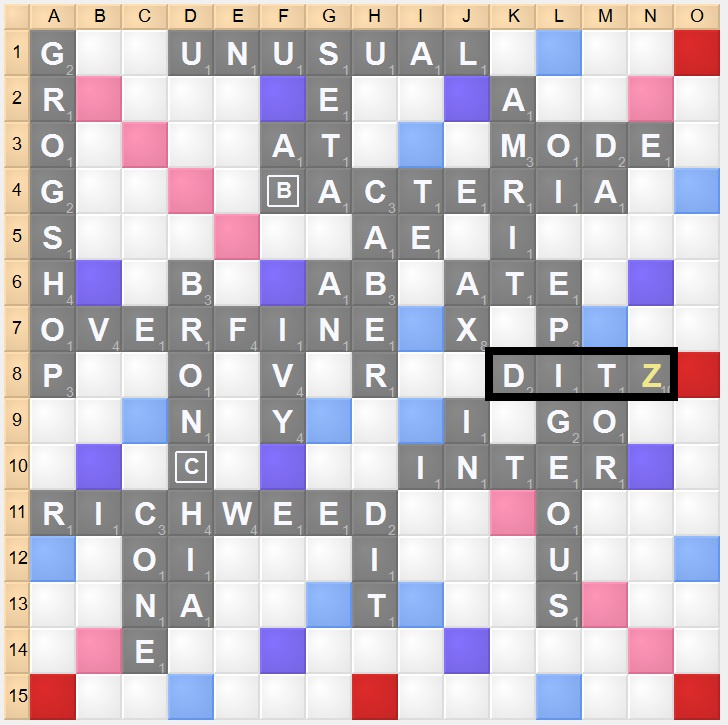}
\ \ \ \ 
\includegraphics[scale=.45]{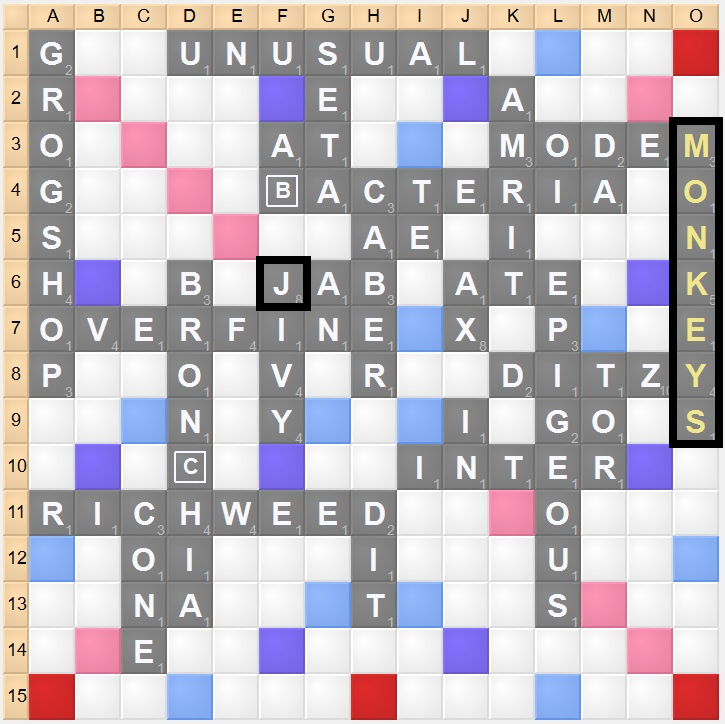}
\end{center}
\end{figure}

The first crucial observation is that if you play
DITZ at 8K (i.e., placing your Z on the 8N square),
then you \emph{set up} a big play down Column~O.
Specifically, if you are lucky enough to draw the E from the bag,
and your opponent makes her highest-scoring play of JIVY at~F6,
then you can bingo out with MONKEYS at~O3,
hooking the M onto MODE to form MODEM
and hooking the Y onto DITZ to form DITZY,
leading to a comfortable 569--537 victory.
See \autoref{fig:dream}.

Although a long shot, this plan seems to give you some winning chances.
``But wait,'' you say.  ``If I play DITZ at 8K,
then won't my opponent sense what I am doing?
Instead of making her highest-scoring play, she can block me,
for example by playing an L at N4, forming the words EL and BACTERIAL.
Then even if I draw the E, I won't be able to play MONKEYS, and I will lose.''

\begin{figure}[ht]
\caption{A 508--508 Tie That Might Worry Your Opponent}
\label{fig:donkeystie}
\begin{center}
\includegraphics[scale=.55]{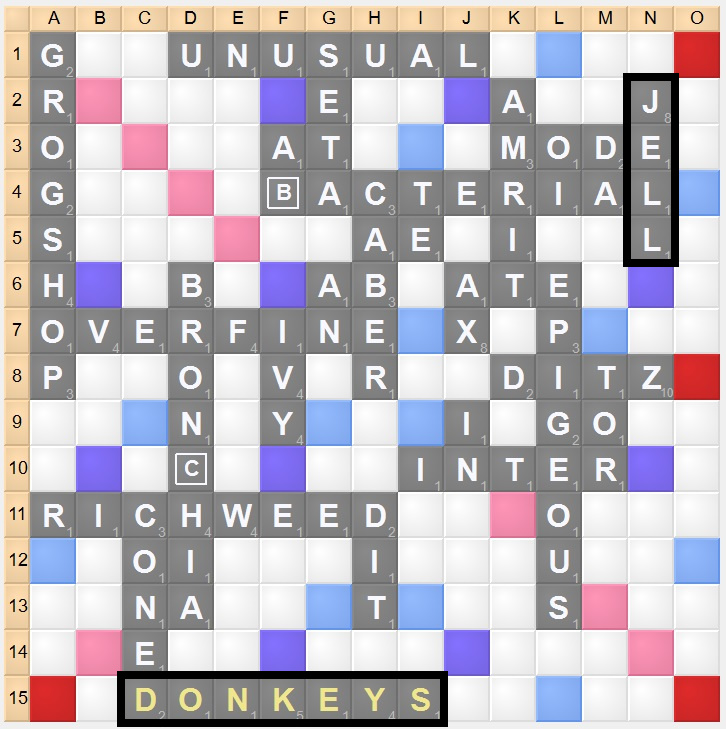}
\end{center}
\end{figure}

We now come to the second crucial observation.
It is important to consider not only the bingos you could draw,
but also the \emph{bingos your opponent worries you might have.}
Specifically, suppose that the last two tiles in the bag are D and~E.
In this scenario, if you play DITZ at 8K,
then the unseen tiles from your opponent's point of view will be
MONKEYSD (i.e., the letters of MONKEYS plus the letter~D).
Your opponent might indeed block MONKEYS;
her highest-scoring blocking play is JELL at N2.
However, she also has to consider the possibility that you are holding DONKEYS,
which plays at 15C, hooking CONE to form CONED.
Is this a serious threat?
If she does the math then she will find that
if you do bingo out with DONKEYS at 15C,
then the result will be a 508--508 tie!
See \autoref{fig:donkeystie}.

Of course, \emph{you} know that this scenario is impossible
(because you're holding the~M and not the~D), but your opponent does not.
Still, at first glance, even this lucky scenario
of drawing the~E and having the~D in the bag seems inadequate.
On your previous turn, you bingoed with GROGSHOP
and drew seven random tiles,
so from your opponent's perspective,
you're equally likely to have drawn the D or the~M at that point,
and MONKEYS poses the bigger threat
(a win as opposed to a tie).
Surely a perfect opponent will block the bigger threat?

\begin{figure}[ht]
\caption{A 511--511 Tie That Could Actually Happen}
\label{fig:monkeystie}
\begin{center}
\includegraphics[scale=.55]{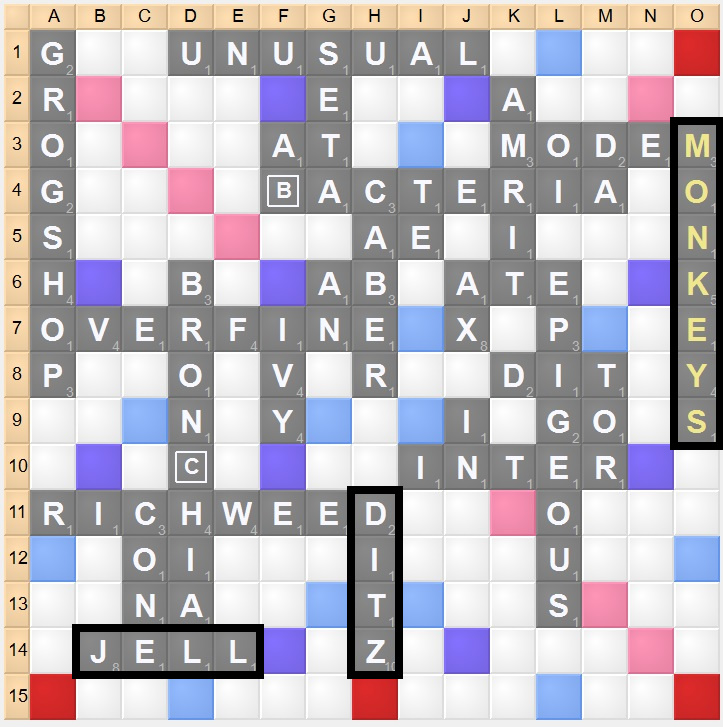}
\end{center}
\end{figure}

Now comes the third crucial observation.
Suppose that the reasoning in the preceding paragraph is correct.
Then what will happen if you play DITZ at H11 instead of at 8K?
Again, if you're lucky and the D and the E are in the bag,
then your opponent will have to worry about both DONKEYS at 15C and MONKEYS at O3.
The difference is that with the Z on the H14 square,
DONKEYS would be the play that hooks DITZ to form DITZY for a big win,
and---as one can check by doing the math---if the opponent plays
her highest-scoring DONKEYS-blocking play of JELL at~14B,
MONKEYS at O3 would lead to a 511--511 tie.
If we believe that the opponent must ``block the bigger threat,''
then she will indeed play JELL at~14B,
giving you free rein to play MONKEYS and tie the game!
See \autoref{fig:monkeystie}.

To summarize, the tentative suggestion is that your best play is DITZ at~H11
(we call this play a \emph{bluff} because
it appears to threaten DONKEYS when you know you cannot play DONKEYS).
The probability is 1/36 that D and~E are in the bag
(because your opponent just bingoed and drew seven random tiles),
and conditional on D and~E being in the bag,
the probability is 1/2 that you will draw the~E.
The argument is that your opponent's best response is to block
DONKEYS by playing JELL at~14B,
letting you tie the game with MONKEYS.
Hence, you have a 1/72 chance of tying and a 71/72 chance of losing.

The above reasoning sounds plausible, but objections can be raised.
If your ``best play'' is to bluff,
then can't we argue that a perfect opponent will realize this,
and \emph{second-guess} you by
responding to your play of DITZ at H11 with JELL at~N2?
Is it really your opponent's ``best play''
to block the bigger threat, if there is an alternative
that performs better against your best play?

To properly analyze the situation,
we must first re-examine the concept of ``best play.''

\section{What Does ``The Best Play'' Even Mean?}

We have been casually using the term ``best play'' as if it is clear
what it means, but in fact it is a subtle concept,
and deserves a closer look.

In the real world, players are fallible, and it could be that
maximizing your practical chances of winning a game involves
exploiting your opponent's fallibility.
For example, in tournament Scrabble, if I play a \emph{phony}
word on the board (i.e., a word that is not listed in the lexicon),
then it is my opponent's responsibility to \emph{challenge} the play 
by asking a referee to check the lexicon.
If my opponent does not challenge, then the phony stays on the board,
and I score points just as if the phony were valid.
If the opponent challenges and the referee determines that
the word is invalid, then I must retract my play, and I lose my turn.
If the opponent challenges and the word turns out to be valid,
then my opponent is penalized (in North America, she loses her turn;
in international play, she does not lose her turn,
but I am awarded 5 extra points).

Given these tournament rules,
I might sometimes intentionally play a word that I know
or suspect is a phony,
betting that my opponent will not dare to challenge it,
especially if my opponent knows that I have superior word knowledge.
Scrabble players refer to such a ploy as \emph{bluffing}.
However, in the present paper, when we talk about the ``best play,''
we assume that both players have perfect word knowledge
and are capable of accurately calculating all possible ways
the game could play out.
With this assumption, ``bluffing'' in the sense of hoping that the opponent has
imperfect word knowledge does not make sense,
so we reserve the term to mean
a valid play that carries a threat that your
opponent does not know you cannot execute.

It will be convenient for the purposes of this paper
to slightly modify the Scrabble scoring system by
awarding $+1$ for a win, $-1$ for a loss, and $0$ for a tie.
This modification has no material effect on the game,
and makes it clearer that Scrabble is a \emph{zero-sum game},
meaning that the sum of the scores of the players is always zero.
In games such as chess, which are not only zero-sum but
which have no hidden information or built-in randomness,
there is a clear notion of ``best play'';
you should make a move that guarantees victory against any opponent,
or if no such move exists, then you should make a move
(if one exists) that guarantees at least a tie against any opponent.

In Scrabble, however, the situation is complicated by
the random draws from the bag
and the \emph{imperfect information}
that the players have about the state of the game.
Until the bag is empty, you are typically unable to
\emph{guarantee} a win or a tie;
usually, the best you can do is to maximize your \emph{probability}
of winning or tying.
Another issue is illustrated by the following example:
suppose you have a choice between two moves, and
\begin{enumerate}
\item Move 1 guarantees that the game will end in a tie, while
\item Move 2 gives you a 30\% chance of winning, a 30\% chance
of tying, and a 40\% chance of losing.
\end{enumerate}
Which move is better?
Near the end of a tournament, you could be in a must-win situation,
where your only path to winning the tournament is to win the current game.
In that case, tying is just as bad as losing, so you will prefer Move~2
to Move~1.
However, during most of the tournament,
you probably want to maximize your \emph{expected score} or \emph{equity},
meaning your probability of winning minus your probability of losing.
In that case, you will prefer Move~1 to Move~2.
A further complication is that, depending on the tournament situation,
players may give some consideration to maximizing spread.

For the rest of our analysis, we ignore spread, and
\emph{we assume that both players play perfectly,
guaranteeing the highest possible equity against a maximally nefarious opponent.}
These assumptions are plausible if, for example,
the winner of the current game takes home all the prize money,
unless the game is a tie, in which case you each take home half,
regardless of spread points.
Under these assumptions, we can apply some mathematical results from
game theory~\cite[Chapters 2 and~6]{KP}.

In particular, an important idea from game theory is that
players may benefit by employing a \emph{mixed} or \emph{randomized} strategy.
Randomization is a familiar idea in games such as
rock-paper-scissors, in which it is obviously a bad idea
to adopt a \emph{pure} or \emph{deterministic} strategy
of always making the same play,
but there is no reason why randomization cannot be used in Scrabble as well.
There is a key theorem about two-player zero-sum games,
called the \emph{minimax theorem},
according to which there must exist a mixed strategy for Player~1
(call it Strategy~A) and a mixed strategy for Player~2
(call it Strategy~B) that together form an \emph{equilibrium}.
This means that if Player~1 adopts Strategy~A,
then Player~2 cannot profit by deviating from Strategy~B,
and similarly, if Player~2 adopts Strategy~B,
then Player~1 cannot profit by deviating from Strategy~A.

In zero-sum games, ``best play'' is commonly defined to mean ``equilibrium strategy.''
The point is that \emph{nonequilibrium} strategies
violate our intuition about what ``best play'' should mean;
if I can profit by unilaterally deviating,
then it does not seem right to say that I am making my best play.
Moreover, the minimax theorem tells us that in zero-sum games,
the equity of an equilibrium is uniquely determined.\footnote{The
equilibrium strategies themselves are not necessarily unique.
In our puzzle, if the opponent second-guesses,
then she just needs to make sure you cannot win or tie with
the rack that she is betting you have,
and it turns out that she has more than one strategy for ensuring that.
The theorem is just that all equilibria yield the same equity.}
(In non-zero-sum games, the situation can be more complicated;
multiple equilibria may exist, with different equities,
and there may be no clear way to even define the term ``best play.'')

\section{What Is a ``Strategy''?}
\label{sec:strategy}

There is one further subtlety we need to address before we are ready
to present the answer to the puzzle.
We have been using the term \emph{strategy} without defining it precisely.
Strictly speaking, to specify a strategy,
we need to specify what both players should do
in \emph{every possible scenario} of the game.
One way to describe the ``game'' we are playing is as follows:
The board position in \autoref{fig:puzzle} is set up,
and then a referee randomly assigns seven of the remaining tiles
to you and seven to your opponent,
putting the remaining two tiles in the bag.
Then the game proceeds according to the rules of Scrabble,
with you playing first, trailing 344 to~476.
Unfortunately, the number of scenarios in this game is enormous,
and specifying what both players should do in every possible
scenario would be a gargantuan task.
However, since we only want the solution to \autoref{fig:puzzle},
we need only specify what to do in a limited number of scenarios.

At first glance, it might seem that
we need only analyze the scenarios in which the referee
assigns you the rack MKNOSYZ, but matters are not quite so simple,
because the game is one of imperfect information.

As we discussed earlier,
playing DITZ in one of the two possible locations gives you
a chance to save the game if D and~E are in the bag.
Therefore, we need to analyze what your opponent will do with
a rack of FJLLLQW.

Now comes the subtle point.
When your opponent holds FJLLLQW, she does not know your rack. 
So we must analyze not only what you should do with your actual rack
of MKNOSYZ,
but also what you should do with any of the racks
that your opponent, holding FJLLLQW, might think you had when you played DITZ.
One way to see why this is necessary is to note that,
even though all available tiles are equally likely to be on your rack
\emph{before} you make your play
(because you just bingoed with GROGSHOP and drew seven random tiles),
your opponent could in principle draw inferences about what tiles you hold
\emph{after} seeing your play of DITZ.
To calculate these inferences correctly,
we need to analyze your behavior holding any of those racks.

If your opponent holds FJLLLQW, then from her perspective,
there were 28 possible racks you could have had when you played DITZ.
Some of these can be eliminated from further consideration.
She can safely ignore the possibility that you had both the D and the M,
reasoning that you cannot now have a bingo on your rack
(see Appendix~\ref{subsec:bingos} for more details).
If you had neither the D nor the~M on your rack when you played DITZ,
then you must have been holding EKNOSYZ,
and it turns out that ZONKEYS is valid
(a zonkey is a cross between a zebra and a donkey).
The puzzle has been carefully crafted so that with a rack of EKNOSYZ,
playing ZONKEYS at~N8 immediately for 126 points is better than playing either DITZ.
Your hope is to draw DL
(not such a long shot given that there are three L's unseen),
for a 537--537 tie after your opponent plays JIVY at~F6 and you go out with DEL at~14B.
Given that you in fact played DITZ,
your opponent can infer that you did not have EKNOSYZ.

We are reduced to considering 12 possible racks for you:
six of them are ``D-racks,'' where you hold DZ plus five of the tiles E, K, N, O, S, Y,
and six of them are ``M-racks,'' where you hold MZ plus five of the tiles E, K, N, O, S,~Y.
See \autoref{tab:12racks}.

\begin{table}[ht]
\caption{The 12 Critical Racks}
\label{tab:12racks}
\centering
\begin{tabular}{c|c}
M-racks & D-racks \\
\hline
\strut
MEKNOSZ & DEKNOSZ \\
MEKNOYZ & DEKNOYZ \\
MEKNSYZ & DEKNSYZ \\
MEKOSYZ & DEKOSYZ \\
MENOSYZ & DENOSYZ \\
MKNOSYZ & DKNOSYZ
\end{tabular}
\end{table}

\section{The Equilibrium Strategies}

We are now ready to describe the equilibrium strategies for both players.

\begin{enumerate}
\item If you hold an M-rack (as in the actual puzzle),
then you play DITZ at H11 with probability 2/3,
and you play DITZ at 8K with probability 1/3.
If you hold a D-rack,
then you play DITZ at 8K with probability 2/3,
and you play DITZ at H11 with probability 1/3.
\item If your opponent holds FJLLLQW with DONKEYSM unseen, then
\begin{itemize}
\item if DITZ is played at 8K, then she plays JELL at N2 with probability 2/3,
and plays JELL at 14B with probability 1/3; and
\item if DITZ is played at H11, then she plays JELL at 14B with probability 2/3,
and plays JELL at N2 with probability 1/3.
\end{itemize}
\end{enumerate}

We can informally summarize the above recommendations as follows.
You should \emph{bluff} with probability 2/3, and \emph{set up} with probability 1/3.
(Here, to \emph{set up} means to play DITZ
in the location that will win the game for you
if you are able to play your bingo,
and to \emph{bluff} means to make the other DITZ play.)
Your opponent should \emph{block} with probability 2/3
and \emph{second-guess} with probability 1/3.
Intuitively, the reason for the 2-to-1 ratio
is that a win is twice as valuable as a tie.

Let us first work out what happens if both players follow the
above recommendations.
We assume that, after you play DITZ,
you draw the tile you need, and leave either a D or an~M in
the bag---whichever one you do not already have on your rack---since
otherwise you simply lose the game.
Then there are four ways the game could play out. For example, with a D-rack:
\begin{enumerate}
\item 8K DITZ, N2 JELL,  probability $(2/3) \times (2/3) = 4/9$, you tie.
\item 8K DITZ, 14B JELL,  probability $(2/3) \times (1/3) = 2/9$, you lose.
\item H11 DITZ, 14B JELL,  probability $(1/3) \times (2/3) = 2/9$, you lose.
\item H11 DITZ, N2 JELL,  probability $(1/3) \times (1/3) = 1/9$, you win.
\end{enumerate}
Your overall equity
(again, conditioned on drawing the tile you need,
and leaving the D or~M in the bag)
is therefore $0 - 2/9 - 2/9 + 1/9 = -1/3$.
Your opponent's equity is of course $+1/3$,
the negative of your equity.

Now suppose that you follow the recommended strategy,
but your opponent decides to deviate, hoping to do better than $+1/3$.
One can check that she cannot profit by deviating.
If she decides to always block, then $2/3$ of the time, she will tie
(from one of your D-racks with 8K DITZ, N2 JELL;
or from one of your M-racks with H11 DITZ, 14B JELL),
and $1/3$ of the time, she will win
(from one of your D-racks with H11 DITZ, 14B JELL;
or from one of your M-racks with 8K DITZ, N2 JELL), and her equity will still be $+1/3$.
If she decides to always second-guess, then $2/3$ of the time, she will win,
and $1/3$ of the time, she will lose, again for an equity of $+1/3$.
Moreover, any randomized mixture of blocking and second-guessing
will also give her an equity of $+1/3$.

On the other hand, suppose your opponent follows the recommended strategy,
and you try to improve on your equity of $-1/3$ by deviating.
It turns out---although it is not trivial to check
(again, see Appendix~\ref{app:analysis})---that
for all 12 critical racks,
you cannot do better than making one of the DITZ plays.
Moreover, deviating from the recommended probabilities does not help.
If you always bluff, then $2/3$ of the time, you will tie
(when your opponent blocks), and $1/3$ of the time, you will lose
(when your opponent second-guesses), for an equity of $-1/3$.
Or if you always set up, then $2/3$ of the time, you will lose,
and $1/3$ of the time, you will win, again for an equity of $-1/3$.
And again, any randomized mixture of bluffing and setting up
will also yield an equity of $-1/3$.

By contrast, the proposed strategies in \autoref{sec:preliminary}
(where you always bluff and the opponent always blocks)
do not form an equilibrium.
We noted this fact already; if you always bluff but your opponent deviates
by always second-guessing, then she profits.
If we think about it for a moment,
it makes intuitive sense that any sound bluffing strategy should be randomized;
the point of bluffing is to exploit your opponent's uncertainty,
so bluffing 100\% of the time would defeat the purpose.

To sum up, with our carefully chosen definition of ``best play,''
we have now justified our initial claim that
your best play in the given position is to bluff with probability 2/3,
and set up with probability 1/3.
Your opponent's best response is to block with probability 2/3
and second-guess with probability~1/3.

\section{Concluding Remarks}

One potentially confusing point is that in \autoref{fig:puzzle},
if your opponent ignores our recommended equilibrium strategy
and plays JELL at N2 regardless of where you play DITZ,
then your opponent will always win. Doesn't this
contradict our claim that your opponent cannot improve on the
equilibrium strategy by unilaterally deviating?
The answer is no; recall how we defined the ``game'' being played
in \autoref{sec:strategy} as starting with the board position
in \autoref{fig:puzzle}, with a referee distributing tiles
randomly to you and your opponent.
We claim only that there is no strategy \emph{for this game as a whole}
that does better than the equilibrium strategy.
It is possible for a deviant strategy to do better
\emph{in certain situations that arise}---but only at the cost
of doing worse in other situations that could arise.

The equilibrium strategies that we have described are not unique.
Our recommendation of a (2/3, 1/3) bluff/setup probability from
every critical rack is the most symmetric,
but another possibility is for you to always bluff from four of your M-racks
and always set up from the other two M-racks
(and similarly for your D-racks).

That one's best play in a zero-sum game with imperfect information can involve
bluffing with some probability strictly between 0 and~1 is of course
not a new idea. Bluffing in poker is a fundamental part of the strategy,
and the mathematical soundness of bluffing was long ago demonstrated in simplified
versions of poker;
see for example ``A Simple Bluff'' in Peter Winkler's book,
\textit{Mathematical Mind-Benders} \cite[page~77]{Win}.
However, to our knowledge, it has never been explicitly demonstrated
before that probabilistic bluffing in Scrabble
can be the best strategy, even against a perfect opponent.

Expert Scrabble players routinely make probabilistic inferences
about the unseen tiles on the opponent's rack.
They may reason, for example, that
``my opponent surely did not have an S when making that play,
because if she did, then she would have bingoed
instead of playing just six tiles.''
Indeed, some experts, by making multiple inferences from
the opponent's plays,
can substantially narrow the range of possible racks that
a strong opponent could have.

Now, whenever inferences are being made,
the possibility of bluffing---i.e., making plays that
induce the opponent to make incorrect inferences---naturally suggests itself.
Our puzzle demonstrates that there do exist situations in Scrabble
where bluffing in this sense is the equilibrium play.
Still, the significance of our puzzle should not be exaggerated.
We have carefully contrived a position in which the circumstances are just right,
and it seems likely that this type of ``inference bluffing''
is not to be recommended in the majority of practical Scrabble positions.

On the other hand, there are some implications for programmers of
Scrabble bots. Some Scrabble bots are capable of perfectly analyzing
\emph{endgame} positions, meaning positions with no tiles left in the bag.
It is tempting to think that perfectly analyzing \emph{pre-endgame}
positions with just one or two tiles in the bag should not be much harder.
However, our puzzle demonstrates that even when there are only one or
two tiles in the bag, mixed strategies can come into play,
and to our knowledge, no existing Scrabble bot understands mixed strategies.
Moreover, a truly thorough analysis of a pre-endgame position
would require assigning a probability of being in the bag
(as opposed to being on the opponent's rack) to every unseen tile,
and these probabilities can be affected by inferences from
previous plays of the game.
In our puzzle, we sidestepped this issue by having both players
bingo just before the pre-endgame, but of course in a real game,
that will rarely be the case.

\appendix
\section{Appendix: More Scrabble Analysis}
\label{app:analysis}

Ideally, we would like to prove that our puzzle is sound,
if necessary with the help of exhaustive computer analysis.
Unfortunately, while the Scrabble analysis engines that we know about,
namely Quackle and Macondo, have been helpful,
they are unable to fully verify soundness.
In this appendix, we sketch the main lines of analysis;
for more details, see \url{https://nackbg.com/scrabble} or
\begin{center}
\url{https://timothychow.net/EquilibriumPuzzle_extended.pdf}
\end{center}

As mentioned in the main part of this paper,
one needs to analyze how your opponent should respond to DITZ
while holding FJLLLQW.
There are 28 racks that you could have had when you played DITZ.
Under the assumption that your opponent holds FJLLLQW, we claim:
\begin{enumerate}
\item if you had neither the D nor the M (i.e., you had ZONKEYS)
then you would have played ZONKEYS rather than DITZ;
\item if you had both the D and the M (15 racks) then your opponent
can guarantee a win by playing JELL at either N2 or 14B; and
\item if you had one of the 12 critical racks, then unilaterally
deviating from our proposed equilibrium strategy would not
have gained you anything.
\end{enumerate}
The first claim we already explained in the main part of the paper.
The second claim is not trivial, but is aided by the
fact that you cannot have a bingo after playing DITZ
(see Appendix~\ref{subsec:bingos} below).
The third claim is the most difficult to verify,
and will be examined more closely below,
but the main idea is that
your opponent can foil your alternative strategies in one of two ways:
she can either \emph{block} your big play,
or (for example if you have two threats that cannot both be blocked)
she can \emph{outrun} you by allowing you to execute your threat but
outscoring you anyway, typically with JIVY at F6, and/or QI at I9.

\subsection{Game Transcript}
\label{subsec:gcg}

We begin by providing a transcript of the game,
in a standard format known as .gcg format.
This proves that the board position can indeed be reached
in a legal Scrabble game, with no phonies played,
with the score being 344 to 476 as we stated.
Each row of the .gcg file specifies which player is making the play,
the player's rack (or part of the player's rack),
the location of the play, the word played
(with dots indicating tiles already on the board,
and lowercase letters indicating blanks),
the points scored by the play, and the player's running total.
The validity of all the words can be confirmed at
\texttt{scrabblecheck.com}.

\begin{verbatim}
#character-encoding UTF-8
#player1 Player_1 Player 1
#player2 Player_2 Player 2
>Player_1: ABCER H4 CABER +24 24
>Player_2: ?AAEIRT 4F bA.TERIA +68 68
>Player_1: AIMT K2 AM.IT +14 38
>Player_2: E 6K .E +2 70
>Player_1: EGIOPSU L6 .PIGEOUS +74 112
>Player_2: O 3K .O +6 76
>Player_1: OT M8 TO +9 121
>Player_2: A 6G A. +4 80
>Player_1: EFINORV 7A OVERFIN. +69 190
>Player_2: D 8K D.. +4 84
>Player_1: ?ABHINO D6 B.ONcHIA +76 266
>Player_2: AX J6 AX +16 100
>Player_1: CDEEIRW 11A RIC.WEED +84 350
>Player_2: NO C11 .ON +16 116
>Player_1: D 3K ..D +18 368
>Player_2: R M8 ..R +5 121
>Player_1: INT 10I INT.. +7 375
>Player_2: AT 3F AT +7 128
>Player_1: VY F7 .VY +9 384
>Player_2: IT H11 .IT +5 133
>Player_1: E 5H .E +4 388
>Player_2: I J9 I. +2 135
>Player_1: E 3K ...E +7 395
>Player_2: E C11 ...E +6 141
>Player_1: E G2 E.. +3 398
>Player_2: GGHOPRS A1 GROGSH.P +203 344
>Player_1: ALNSUUU 1D UNUSUAL +78 476
\end{verbatim}

\subsection{Bingos You Could Draw}
\label{subsec:bingos}

We start by assuming that you need to bingo to win.
See Appendix~\ref{subsec:needtobingo} below for a
discussion of the possibility of winning without bingoing.

Using a computer, we can search the lexicon for bingos
that you could hope to play, starting with 
any of the 12 critical racks.
Besides DONKEYS, MONKEYS, and ZONKEYS, the list includes
the following seven-letter words:
\begin{equation*}
\begin{array}{llll}
\text{DOLMENS} &\text{DONZELS} &\text{ENFOLDS} &\text{FONDLES} \\
\text{KNOLLED} &\text{MENFOLK} &\text{ZEDONKS}
\end{array}
\end{equation*}
As for eight-letter words, from the rack MENOSYZ, you could play DITZY at 8K and draw LL,
threatening SOLEMNLY at O1.
Or from the rack DEKOSYZ, you could play ZEK at G10 and draw LM,
and if your opponent were to float an L by playing MODEL at 3K,
then you could play SELDOMLY at~O1.

There are other bingos that can be formed from the available letters.
The seven-letter bingos are FELLOWS, MELLOWS, SELFDOM, SWOLLEN, and YELLOWS.
The longer bingos are FELLOWLY, MELLOWLY, MENFOLKS, and
MINDFULLY (at 12G, through the existing tiles I and~U).
However, these bingos
cannot actually be played starting from any of the 12 critical racks.

To confirm a claim we made earlier,
we point out that all the bingos containing both D and~M also contain an~L;
therefore, if you have a rack containing DMZ and four of the
letters E, K, N, O, S, Y, then you cannot draw into a bingo by playing DITZ
if your opponent holds FJLLLQW.

\subsection{Bingoing Immediately}
\label{subsec:zedonks}

Since you need to bingo to win, one option is to bingo immediately.
Only one of the 12 critical racks admits this option, namely DEKNOSZ,
which allows you to play ZEDONKS immediately at 15A
(but not at O1, since ``MODED'' is not valid).
This play looks promising; it almost levels the score (468--476).
However, it turns out that there is no way you can tie or win.

If you draw the~J, then your opponent will play ABLATE at 6G,
blocking you from playing JIVY at F6.
You are now ``J-stuck'' (meaning that you have nowhere to play your~J)
and your opponent, who has multiple spots to play her~Q,
will coast to an easy win.

If you do not draw the~J, then your opponent can play JIVY at F6
(or, if she has the right tiles, JELLY at N2)
and outrun you no matter what you draw.
For example, if you draw LL then you can go out with ELL at~N3,
losing 533--537.

\subsection{Two-Tile Plays}
\label{subsec:twotiles}

Since there are only two tiles in the bag, the only way you can bingo
(if you do not bingo immediately) is to play either two tiles or one tile.
In this section, we consider two-tile plays.
A general problem with a two-tile play is that it empties the bag,
giving your opponent full knowledge of your rack,
which usually simplifies her task of countering your threats.
Often she has a choice of either blocking or outrunning.

For example, let us consider the original rack MKNOSYZ.
You could play MY at 2N hoping to draw into ZEDONKS,
but ZEDONKS plays only at 15A,
so your opponent can either outrun with JIVY,
or block the only spot where your bingo plays.

There are some cases in which you have a bingo that plays in two spots,
forcing your opponent to outrun. For example,
with the rack DEKNOSZ, you could play ZEK at G10 and draw into
DOLMENS, which plays at both O1 and~15C.

On the other hand, if you play DITZY at either 8K or H11 for 54 points,
then outrunning may not be possible.
However, after DITZY, any bingo you have can be blocked.
The Y of DITZY will interfere with bingos that
might otherwise play along that edge of the board,
except for SOLEMNLY, which \emph{requires} the Y of DITZY.

\subsection{One-Tile Plays}
\label{subsec:onetile}

Apart from DONKEYS and MONKEYS,
the bingos that you could draw into with a one-tile play from
one of the 12 critical racks are DONZELS, ZONKEYS, and ZEDONKS.
Again, we give some representative lines rather than
a comprehensive analysis.

The rack DEKNOSZ allows you to try for DONZELS---which plays
at both 15C and 14E---by playing KI at I9 and drawing an~L.
If the M is in the bag, then you will also appear to threaten
DOLMENS at O1, but your opponent can outrun with JIVY or JELLY.
If the J is in the bag then DONZELS is your only threat,
and she can block both spots with 14C~EL.

You can try to draw into ZONKEYS from an M-rack by playing MI at 9I
for 14 points (playing off a D from a D-rack leads to similar considerations).
\begin{enumerate}
\item Suppose you hold MKNOSYZ, play MI, draw the E, and leave the D in the bag.
Then your opponent cannot block both ZONKEYS at N8 and DONKEYS at 15C.
However, she can play JIVY for 61 points,
and after you bingo out with ZONKEYS for 126 points
plus 46 points for your opponent's rack of DFLLLQW,
you still lose 530--537.
\item Similarly, suppose you hold MEKNOSZ, play MI, draw the D,
and leave the Y in the bag.
Then your opponent cannot block both ZONKEYS at N8 and ZEDONKS at 15A,
but again JIVY outruns both; ZEDONKS loses 532--537.
\item Suppose you hold MEKOSYZ, play MI, draw the J, and leave the N in the bag.
You have failed to draw ZONKEYS,
but ZONKEYS is still a threat from your opponent's perspective,
and if she is careless about how she blocks it,
then there is some risk that you could use your J and Z to score
enough points to win.
Fortunately for her, she can win by playing WILLFUL at 12G;
you can play ZEK at 14H, preparing to go out with JOYS at 15F,
but your opponent can play QI at B10 and win by 23 points.
Or if instead of ZEK, you try OK at 10A to block QI,
then your opponent can still win by playing QAT at 3E.
\end{enumerate}

\subsection{Non-Bingo Variations}
\label{subsec:needtobingo}

In \autoref{sec:preliminary}, we claimed that with a deficit of 132 points,
you need to bingo in order to avert a loss.
Exhaustively checking all non-bingo variations
is not realistically possible by hand,
and as we said earlier, Quackle and Macondo are not equipped for the task either.
However, we have done enough manual analysis to be
highly confident that there is no way to
save the game without bingoing.
In this section, we describe some ``near-misses''
to illustrate the types of plays that need to be examined closely.

If you have a rack of MEKOSYZ and you play DITZ at H11,
then you could draw a J and leave N in the bag.
Your opponent sees the possible threat of MONKEYS
and can block with DWELL at N1.
You have no bingo and will lose,
but if you play SMOKEY at 15C,
then your opponent's best response is ABLATE at 6G,
allowing you to go out with JO at E14 for a 3-point loss.
(We should point out that instead of blocking with DWELL,
your opponent could outrun with QI at I9,
winning by a larger margin in the non-bingo variations,
at the cost of barely outrunning MONKEYS.)

If you have a rack of MKNOSYZ, then you can play ZONK at N8 for 45 points.
Your opponent counters with JIVY and then you go out
with EMYDS at O6, losing 531--537.

For a different type of near-miss, the rack MENOSYZ would yield
winning chances if ``ZON'' at N8 were valid.
You would hope to draw into MYSELF, which plays at
O1, O7, and~15A.

\subsection{Miscellaneous Final Comments}
\label{subsec:misc}

Any of the 12 critical racks could be used to pose a valid puzzle.
We chose MKNOSYZ because it minimizes the number of bingos you could try for.
For an audience of Scrabble experts, we prefer the rack DEKNOSZ, 
in part because of the tempting possibility of playing ZEDONKS immediately,
and in part because it is surprising to set up a hook for a Y that
you do not possess.

ZEDONKS and ZONKEYS were added to NWL2023 and were not valid
in earlier versions of the lexicon.
Suppose ZONKEYS were not valid and
the only way to save the game with a rack of ZONKEYS
were to play DITZ.
You would be equally likely to draw DONKEYS or MONKEYS
after playing off the~Z,
so you would have no reason to prefer one DITZ over the other.
This change would lead to slightly messier equilibrium probabilities,
but if one is not bothered by that,
then one might be able to compose a variation of the puzzle
that has fewer bingo possibilities and that is therefore easier to check
for soundness. For example, instead of DONKEYS/MONKEYS
one could try FUSTILY/JUSTIFY and an opponent rack of CUUUVVX.

\begin{figure}[ht]
\caption{\texttt{https://www.cross-tables.com/annotated.php?u=27771\char35 5}}
\label{fig:sokol}
\begin{center}
\includegraphics[scale=.5]{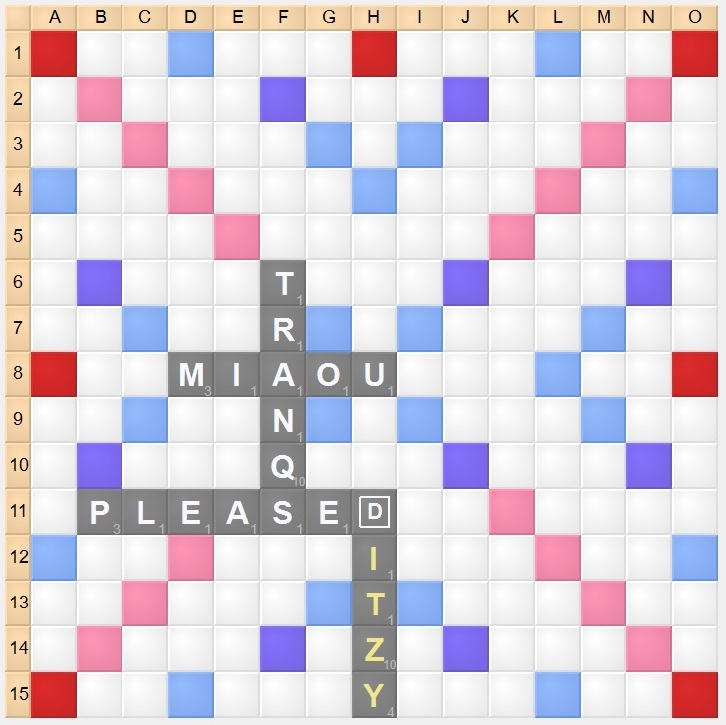}
\end{center}
\end{figure}

After we composed this puzzle, Jerry Lerman informed us of a
tournament game that he played against Josh Sokol in which
Sokol, holding AIILTYZ, actually played DITZY at H11! See \autoref{fig:sokol}.
Simulations by Quackle and Macondo suggest that maybe
Sokol should have played DITZ,
setting up his~Y instead of cashing it immediately.
We thank Jerry, as well as Charlie Carroll and Joe Edley,
for valuable comments on an early version of this work.

The equilibrium we described almost works in the other
major English-language lexicon currently in use, namely CSW2024.
However, there is a subtle problem;
if you started with the rack DENOSYZ,
then a better play than DITZ is ZEN at N8.
GOE\# is valid in CSW2024
(the \# indicates that the word is valid in CSW but not NWL),
and if you draw JL, then you win.
If your opponent tries QI at I9,
then you play ODYLS at O6;
your opponent's best try is WELK\# at 14B,
but then you play JIVY at F6 and win by 1~point.
Or if instead your opponent tries GOEL\# at 9L,
then the optimal sequence is
JIVY at F6, WELK\# at 14B, ODYLS at O6,
QI at I9, and EL at~N3, and again you win by 1~point.
Again, the equilibrium probabilities for CSW2024
are therefore a bit messier than 1/3 and~2/3.

Not long after the present paper was mostly written,
we devised a version of the puzzle that works
with both the NWL2023 and CSW2024 lexicons.
This version was published in NASPA News
as the ``Equilibrium Puzzle'' in late 2025;
the interested reader can find it at
\url{https://timothychow.net/misc/NASPA-2025-11-06.html}

The appearance of the word UNUSUAL in \autoref{fig:puzzle} is
something of an inside joke; it is an allusion to
a puzzle by the first author that also features the word UNUSUAL. 
See \url{https://nackbg.com/scrabble} or
\begin{center}
\url{https://timothychow.net/UnusualPuzzle_extended.pdf}
\end{center}

\bigskip\hrule\bigskip
\noindent
\textbf{Email:} \texttt{tchow@alum.mit.edu, nack2000@sbcglobal.net}


\begin{thebibliography}{1}
\bibitem{KP} Anna R. Karlin and Yuval Peres,
\textit{Game Theory, Alive}, American Mathematical Society, 2017.
\bibitem{Win} Peter Winkler,
\textit{Mathematical Mind-Benders}, A K Peters, 2007.
\end{thebibliography}
\end{document}